\documentclass[10pt]{article}
\usepackage{graphicx}
\usepackage{amsfonts}
\usepackage{amsmath}
\usepackage{amssymb}
\usepackage{eufrak}
\usepackage{color}
\usepackage{algorithm}
\usepackage{algorithmic}
\usepackage{rotating}
\usepackage{textcomp}
\usepackage{hyperref}% HAS to be the last one defined
\newtheorem{thm}{Theorem}[section]
\newtheorem{cor}{Corollary}[section]

\newcommand{\R}{\mathbb R}

\newcommand{\G}{\mathcal{G}}

\newcommand{\OO}{\mathcal{O}}

\renewcommand{\L}{\mathcal{L}}

\newcommand{\la}{\left\langle}
\newcommand{\ra}{\right\rangle}
\newcommand{\lno}{\left\|}
\newcommand{\rno}{\right\|}
\newcommand{\de}{d_{\Delta_n}}
\newcommand{\boundReturns}{L}

\newcommand{\dualElement}{g}

\begin{document}
\title{Hedge algorithm and Dual Averaging schemes}
\author{Michel Baes\footnote{Institute for Operations Research, ETH Z\"urich, R\"amistrasse 101, 8092 Z\"urich, Switzerland,\newline michel.baes@ifor.math.ethz.ch.} , Michael B\"urgisser\footnote{Institute for Operations Research, ETH Z\"urich, R\"amistrasse 101, 8092 Z\"urich, Switzerland,\newline michael.buergisser@ifor.math.ethz.ch.}}

\maketitle

\begin{abstract}
We show that the Hedge algorithm, a method that is widely used in Machine Learning, can be interpreted as a particular instance of Dual Averaging schemes, which have recently been introduced by Nesterov for regret minimization. Based on this interpretation, we establish three alternative methods of the Hedge algorithm: one in the form of the original method, but with optimal parameters, one that requires less a priori information, and one that is better adapted to the
context of the Hedge algorithm. All our modified methods have convergence results that are better or at least as good as the performance guarantees of the vanilla method. In numerical experiments, our methods significantly outperform the original scheme.
\end{abstract}

\section{Introduction}

The Hedge algorithm was introduced by Freund and Schapire
\cite{Freund_Hedge} and encompasses many well-known schemes in
Machine Learning. For instance, as Freund and Schapire showed, this
method is related to the now widely used AdaBoost algorithm
\cite{Freund_Hedge}. The Hedge algorithm can be used to solve the
following online allocation problem. We want to invest an amount of
money in a portfolio consisting of different assets at the stock
market. After each time step, we can modify the current composition
of our portfolio. The Hedge algorithm defines an update strategy for
our portfolio, such that the average performance that we achieve is
not much worse than the average performance of the most favorable
investment product. The portfolio update rule is based on the current
loss (or gain) that is associated with every investment product.

In this paper, we propose an alternative viewpoint on the Hedge
algorithm, using methods that have recently been introduced in Convex
Optimization. It is well-known that the Hedge algorithm can be
interpreted as a Mirror-Descent scheme \cite{Nemirovski_Yudin_book}
with an entropy-type prox-function; see for instance Chapter 11 in
\cite{Cesa-Bianchi}. However, this interpretation has two drawbacks.
First, Mirror-Descent schemes require the definition of a convex and
closed objective function. In this setting, the current loss of the
investment products corresponds to a subgradient of this objective
function. In particular, we explicitly rule out the possibility of a
dynamic objective function with this approach. However, modeling the
performance of a portfolio with a static objective function, even
when we allow random losses, is at best questionable. As the last
financial crisis has shown, significant sudden changes in the
performance of an investment product can appear, that are more
appropriately modeled with a dynamic objective function. Second, in
order to ensure convergence, Mirror-Descent schemes need to consider
subgradients with more weight the earlier they appear. However,
common sense dictates that recent losses contain more relevant
information on the future development of the stock market than losses
occurred years ago. In this paper, we interpret the Hedge algorithm
as a Dual Averaging scheme \cite{nesterov:PDS}.

Dual Averaging schemes are the natural extension of Mirror-Descent
methods and get rid of both deficiencies we pointed out above at the
same time. When applied to our context, Dual Averaging schemes do not
make any assumptions on the construction of the losses. For instance,
they can be chosen in adversarial way with respect to our current
portfolio, they can be randomly generated, or \--- which reflects
some of the latest events at the stock market more accurately \---
their construction rule may dynamically change. Moreover, in Dual
Averaging schemes, we can give more weight to the latest losses,
which allows to react much faster to significant changes in the
market behavior.

Based on this alternative interpretation of the Hedge algorithm, we
give three modifications of the Hedge algorithm, namely the Optimal
Hedge algorithm, the Optimal Time-Independent Hedge algorithm, and
the Optimal Aggressive Hedge algorithm. All these methods have
convergence results that are better or at least as good as the
convergence guarantee for the vanilla Hedge algorithm. The Optimal
Hedge algorithm has the same form as the original Hedge algorithm,
except that all method parameters are chosen in an optimal way. The
Optimal Time-Independent Hedge algorithm requires less a priori
information than the Optimal Hedge algorithm. Finally, the Optimal
Aggressive Hedge algorithm considers losses as more relevant the
later they appear. Numerical results show that all our alternative
methods perform better than the vanilla Hedge algorithm. More
interestingly, using the Optimal Aggressive Hedge algorithm, we end
up with an average benefit that is even better than the profit of the
single most favorable investment product, provided that the losses
incur shocks reverting the performance of assets. This effect would
not have been possible with a static objective function.

This paper is organized as follows. In Sections \ref{sec:DAS} and \ref{sec:Hedge}, we review Dual Averaging schemes and the original Hedge algorithm. We show in Section \ref{sec:Inter} that the Hedge algorithm is a Dual Averaging scheme and suggest several alternative methods based on this interpretation. We conclude this paper with some numerical results in Section \ref{sec:Num}.

\section{Dual Averaging methods}
\label{sec:DAS}

We give a brief review of Dual Averaging schemes, which were
introduced by Nesterov in \cite{nesterov:PDS}.

Let $Q\subset\R^n$ be a closed and convex set. We assume that we have
at our disposal an oracle $\G$, which returns a vector
$\dualElement=\G(x)\in\R^n$ for input $x\in Q$. We interpret the
oracle output $\dualElement=\G(x)$ as a loss vector that is
associated to $x$. The corresponding loss is defined as $\la
\dualElement,x\ra$, where $\la\cdot,\cdot\ra$ denotes the standard
dot product in $\R^n$. Assume now that we repeat this process. That
is, for $t\in\mathbb{N}$, we choose an element $x_t\in Q$, call the
oracle $\G$ with input $x_t$, observe the loss vector
$\dualElement_t=\G(x_t)\in\R^n$, and update our choice of the element
$x_{t+1}\in Q$. After $T$ rounds, we obtain a total averaged loss of
\[\L_T:=\frac{1}{\sum_{k=0}^{T-1}\lambda_k}\sum_{t=0}^{T-1}\lambda_t\la \dualElement_t,x_t\ra,\]
where the numbers $\lambda_0,\ldots,\lambda_{T-1}>0$ can be seen as a tool to weight the losses according to their appearance. We can compare $\L_T$ to the averaged loss $\sum_{t=0}^{T-1}\lambda_t\la \dualElement_t,\bar{x}\ra/\sum_{k=0}^{T-1}\lambda_k$, where $\bar{x}$ corresponds to an element in $Q$ that turns out to be optimal in hindsight. The deviation of this two quantities is called \textit{averaged regret} and denoted by $\mathcal{R}_T$:
\[\mathcal{R}_T:=\frac{1}{\sum_{k=0}^{T-1}\lambda_k}\left(\sum_{t=0}^{T-1}\lambda_t \la \dualElement_t, x_t\ra-\min_{x\in Q} \left\{\sum_{t=0}^{T-1}\lambda_t \la \dualElement_t, x \ra\right\}\right)=\frac{1}{\sum_{k=0}^{T-1}\lambda_k}\max_{x\in Q}\left\{ \sum_{t=0}^{T-1}\lambda_t\la \dualElement_t,x_t- x \ra\right\}.\]
If the oracle $\G$ is associated to a convex optimization problem of the form $\min_{x\in Q}f(x)$, that is, the oracle return correspond to subgradients of $f$, the averaged regret $\mathcal{R}_T$ gives us an upper bound on the optimality gap $\min_{0\leq t\leq T-1}f(x_t)-\min_{x\in Q}f(x)$.

Naturally, the following question arises: is there a strategy to update the elements $x_0,\ldots,x_{T-1}$ such that the averaged regret $\mathcal{R}_T$ is bounded from above by a quantity that converges to zero when $T$ goes to infinity? Nesterov's Dual Averaging schemes \cite{nesterov:PDS} can be used to define such update strategies.

We equip $\R^n$ with a norm $\lno\cdot\rno$, not necessarily the norm associated with $\la\cdot,\cdot\ra$, and denote by $\lno\cdot\rno_*$ the corresponding dual norm. Nesterov's Dual Averaging methods require a prox-function $d:Q\rightarrow\R$, that is, a function that is continuous and strongly convex modulus $\sigma>0$ with respect to $\lno\cdot\rno$ on $Q$. We set $x_0=\arg\min_{x\in Q}d(x)$. Without loss of generality, we may assume that $\sigma=1$ and that $d$ vanishes at $x_0$. The algorithm accumulates all the loss vectors in a dual variable $s_{t+1}$, that is, $s_{t+1}=-\sum_{k=0}^{t}\lambda_k\dualElement_k$ for any $t=0,\ldots,T-1$. In order to define $x_{t+1}$, the dual variable $s_{t+1}$ is then projected back on the set $Q$ using the \textit{parametrized mirror-operator}
\[\pi_{Q,\beta_{t+1}}:\R^n\rightarrow Q: s\mapsto \arg\max_{x\in Q}\left\{\la s,x-x_0\ra-\beta_{t+1} d(x)\right\},\]
where $\beta_{t+1}>0$ is some projection parameter. We assume that $d$ is chosen in such a way that the above optimization problem is easily solvable. The resulting scheme looks as follows.

\begin{algorithm}[h]
  \caption{Dual Averaging methods \cite{nesterov:PDS}}
\label{alg:DAS}
  \begin{algorithmic}[1]
\STATE Set $s_0=0$ and $x_0=\arg\min_{x\in Q} d(x)$.
\STATE Choose positive weights $\{\lambda_t\}_{t\geq 0}$
and a non-decreasing sequence $\{\beta_t\}_{t\geq 0}$ of positive projection
parameters.
\FOR{$t=0,\ldots,T-1$}
\STATE Call the oracle $\mathcal{G}$ to get a loss vector $\dualElement_t=\mathcal{G}(x_t)\in\R^n$.\\
\STATE Set $s_{t+1}=s_t-\lambda_t\dualElement_t$.\\
\STATE Compute $x_{t+1}:=\pi_{Q,\beta_{t+1}}(s_{t+1})$.\\
\ENDFOR
\end{algorithmic}
\end{algorithm}

Nesterov proved the following result for this method.

\begin{thm}[First part of Theorem 1 in
\cite{nesterov:PDS}]\label{thm:ConvergenceDAS}
For any $D\geq 0$, we have:
\begin{equation}
\label{eq:ConvergenceDAS}
\frac{1}{\sum_{k=0}^{T-1}\lambda_k}\max_{x\in Q}\left\{\sum_{t=0}^{T-1}\lambda_t\left\langle \dualElement_t,x_t-x
\right\rangle:d(x)\leq D\right\}\leq \frac{1}{\sum_{k=0}^{T-1}\lambda_k}\left(\beta_{T}D +
\frac{1}{2}\sum_{t=0}^{T-1}\frac{\lambda_t^2}{\beta_{t}}||\dualElement_t||^2_*\right).
\end{equation}
\end{thm}

Let us assume that the oracle returns are uniformly bounded, that is, there exist a constant $\boundReturns$ such that $\lno\dualElement_t\rno_*\leq\boundReturns$ for any $t=0,\ldots,T-1$. The above theorem motivates several ways to choose the weights $\lambda_t$ and the projection parameters $\beta_t$. For instance, we can set $\beta_t=1$ for any $t$ and choose  constant weights $\lambda_t =\lambda^*$ in such a way that the right-hand side in (\ref{eq:ConvergenceDAS}) is minimized. That is, provided that $T$ is fixed in advance, we set $\lambda^*=(1/\boundReturns)\sqrt{2D/T}$, for which the right-hand side in (\ref{eq:ConvergenceDAS}) becomes $\boundReturns\sqrt{2D/T}$. %Note that the optimal constant weight choice depends on $T$.
Moreover, Nesterov observed in \cite{nesterov:PDS} that for $\beta_t=1$ the right-hand side in (\ref{eq:ConvergenceDAS}) converges to zero as long as as $\sum_{t=0}^T\lambda_t$ diverges and $\sum_{t=0}^T\lambda_t^2$ converges when $T$ goes to infinity. The latter condition implies that the weights $\lambda_t$ converge to zero. Selecting the $\beta_t$'s in an appropriate way, we can allow non-decreasing weights $\lambda_t$ while still ensuring that the right-hand side in (\ref{eq:ConvergenceDAS}) converges to zero when $T$ goes to infinity. For instance, as Nesterov suggested in \cite{nesterov:PDS}, we can set
\begin{eqnarray}\label{eq:betas}
\lambda_t=\frac{\sqrt{2D}}{\boundReturns},\qquad \beta_0=1,\qquad\text{and}\qquad \beta_{t+1}=\sum_{k=0}^t\frac{1}{\beta_t}\qquad \forall\ t\geq 0,
\end{eqnarray}
for which the right-hand side in (\ref{eq:ConvergenceDAS}) is still in $\OO\left(\boundReturns\sqrt{D/T}\right)$. The same asymptotic bound can be guaranteed for $\lambda_t=(t+1)^2\sqrt{7D}/\boundReturns$ and $\beta_{t}=t^{2.5}$ for each $t\geq 0$.

\section{The Hedge algorithm}
\label{sec:Hedge}

The Hedge algorithm \cite{Freund_Hedge} is a generic method that encompasses many well-known schemes in Machine Learning. As examples, Multiplicative Weights Update methods are a variation of the Hedge algorithm (see \cite{Arora_MWU} for a survey) and AdaBoost can be related to the Hedge algorithm (see \cite{Freund_Hedge} for more details).

The problem the Hedge algorithm aims at solving can be described as follows. We assume that we want to invest a certain amount of money at the stock market. We have at our disposal a basket of $n$ investment products such as shares, currencies, gold, raw materials, real estates, and so on. Let us denote by $x_{t,i}\geq 0$ the share of our initial amount of money that we invest in product $i$ at time $t$, where $i=1,\ldots,n$ and $t\geq 0$. We always invest all of our money, that is, we assume $\sum_{i=1}^nx_{t,i}=1$ for all $t\geq 0$. At every time step $t\geq 0$, we can evaluate the loss (or gain) $\ell_{t,i}$ corresponding to the investment product $i$, where we assume $\ell_{t,i}\in [-\mu,\rho]$ for every $t\geq 0$ and any $i=1,\ldots,n$. Thus, given our portfolio $x_t$ at time $t$, we suffer a loss of $\la \ell_t,x_t\ra$ at this time step. The Hedge algorithm defines now an update strategy for our portfolio such that the averaged loss $\sum_{t=0}^{T-1}\la \ell_t,x_t\ra/T$ that we face is not much worse than the averaged total loss $\min_{1\leq i\leq n}\sum_{t=0}^{T-1}\ell_{t,i}/T$ of the investment product with the best performance.

The Hedge scheme evaluates the losses through a decreasing score function $U:[-\mu,\rho]\rightarrow (0,1]$. For the sake of brevity, we focus in this paper only on score functions of the form $U(z)=\gamma^{az+b}$, where $\gamma\in (0,1)$, $a>0$, and $b\in\R$ are some parameters whose choices we discuss in detail afterwards. The Hedge algorithm assigns a weight $w_{t,i}$ to every investment product $1\leq i\leq n$ and for every time step $t\geq 0$. The current weight of investment product $i$ depends on its initial weight and on its performance in the past. More concretely, it is defined as $w_{t+1,i}:=w_{t,i}U(\ell_{t,i})$. The portfolio $x_{t+1}$ is then given by the normalization of the weight vector $w_{t+1}$. The full method takes the following form.

\begin{algorithm}[h]
  \caption{Hedge algorithm \cite{Freund_Hedge}}
\label{alg:Hedge}
  \begin{algorithmic}[1]
\STATE Let $\gamma\in (0,1)$, $a>0$, $b\in\R$, and $w_0=(1/n,\ldots,1/n)\in\R^n$.
\FOR{$t=0,\ldots,T-1$}
\STATE Set $x_t=w_t/\sum_{i=1}^nw_{t,i}$.\\
\STATE Observe loss vector $\ell_t$.\\
\STATE Set $w_{t+1,i}=\gamma^{a\ell_{t,i}+b}$ for any $i=1,\ldots,n$.\\
\ENDFOR
\end{algorithmic}
\end{algorithm}

Freund and Schapire studied the convergence behavior of Algorithm \ref{alg:Hedge}. In their paper, they considered the situation where $\mu=0$ and $\rho=1$. The immediate extension of their reasoning to our more general setting yields to the following result.

\begin{thm}[Extension of Theorem 2 in \cite{Freund_Hedge}] \label{thm:convergenceHedge} With $a=1/(\mu+\rho)$ and $b=\mu/(\mu+\rho)$, the sequence $\left(x_t\right)_{t=0}^{T-1}$ generated by Algorithm \ref{alg:Hedge} satisfies
\begin{eqnarray}\label{eq:freund_schapire_intermediate}
\sum_{t=0}^{T-1}\left(\mu+ \la \ell_t,x_t\ra\right)\leq \frac{\mu+\rho}{1-\gamma}-\frac{\ln(\gamma)}{1-\gamma}\min_{1\leq i\leq n}\left(\sum_{t=0}^{T-1}\left(\mu+\ell_{t,i}\right)\right).
\end{eqnarray}
\end{thm}

As mentioned in \cite{Freund_Hedge}, the above theorem can be extended to any decreasing score function $U:[-\mu,\rho]\rightarrow\R$ that complies with the condition
\begin{eqnarray}
 \label{eq:condition_score_function}
\gamma^{\frac{z+\mu}{\mu+\rho}}\leq U(z)\leq 1 - (1-z)\frac{z+\mu}{\mu+\rho}\qquad \forall\ z\in [-\mu,\rho].
\end{eqnarray}

Optimizing the right-hand side of (\ref{eq:freund_schapire_intermediate}) with respect to $\gamma$, that is, setting $\gamma=1/\left(\sqrt{2\ln(n)/T}+1\right)$, we obtain the score function
\begin{eqnarray}\label{eq:score_Freund}
U:[-\mu,\rho]\rightarrow\R:z\mapsto \left(\sqrt{2\ln(n)/T}+1\right)^{-\frac{z+\mu}{\mu+\rho}},
\end{eqnarray}
for which one can prove the following statement using Theorem \ref{thm:convergenceHedge}; see \cite{Freund_Hedge} for more details on the derivation of this score function.

\begin{cor}[Consequence of Lemma 4 in \cite{Freund_Hedge}] With the above score function, we have:
\begin{eqnarray}\label{eq:bound_hedge}
\frac{1}{T}\left(\sum_{t=0}^{T-1}\la \ell_t,x_t\ra-\min_{1\leq i\leq n}\sum_{t=0}^{T-1} \ell_{t,i}\right)\leq (\mu+\rho)\left(\frac{\ln(n)}{T}+\sqrt{\frac{2\ln(n)}{T}}\right).
\end{eqnarray}
\end{cor}

\section{The Hedge algorithm is a Dual Averaging method}
\label{sec:Inter}

We show that we can recast the Hedge algorithm in the framework of Dual Averaging schemes and derive alternative versions of the original method.

We define $Q$ as the $(n-1)$-dimensional standard simplex $\Delta_n=\left\{x\in\R^n:x\geq 0, \sum_{i=1}^nx_i=1\right\}$, so that $Q$ encompasses all possible portfolios. We equip $\R^n$ with the norm $\lno x\rno_1:=\sum_{i=1}^n|x_i|$, for which the corresponding dual norm is of the form $\lno s\rno_\infty:=\max_{1\leq i\leq n}|s_i|$. Moreover, we endow Algorithm \ref{alg:DAS} with the prox-function
\[\de:\Delta_n\rightarrow\R:x\mapsto \ln(n)+\sum_{i=1}^nx_i\ln(x_i).\]
It is well-known that this function complies with our assumptions on prox-functions, that is, it is continuous and strongly convex modulus $1$ with respect to $\lno\cdot\rno_1$ on $\Delta_n$, it attains its center $x_0:=\arg\min_{x\in\Delta_n}\de(x)$ at $(1/n,\ldots,1/n)$, and it vanishes at this point; see for instance \cite{nesterov:PDS} and the references therein. Moreover, we can explicitly write the corresponding parametrized mirror-operator:
\[\pi_{\Delta_n,\beta}(s)=\left(\frac{\exp(s_i/\beta)}{\sum_{j=1}^n\exp(s_j/\beta)}\right)_{i=1}^n\qquad \forall\ s\in\R^n,\ \forall\ \beta >0.\]
Given a loss vector $\ell_t\in\R^n$, that is, $\ell_{t,i}$ corresponds to the loss of investment product $i$ at time $t$, we evaluate this vector trough an affine function $z\mapsto az+b$, where $a>0$ and $b\in\R$. We interpret the resulting vector $g_t:=\left( a\ell_{t,i}+b\right)_{i=1}^n$ as the output $\G(x_t)$ of the oracle $\mathcal{G}$ in Algorithm \ref{alg:DAS} for input vector $x_t\in\R^n$. Note that we do not specify any construction rule for the loss vector $\ell_t$. For instance, they could be chosen randomly or in an adversarial way with respect to the portfolio $x_t$. Algorithm \ref{alg:DAS} takes the following form for our setting, where we express the parametrized mirror-operator $\pi_{\Delta_n,\beta}$ in a form that makes the comparison of the resulting method with the Hedge algorithm rather transparent.

\begin{algorithm}[h]
  \caption{Extended Hedge algorithm}
\label{alg:Extended_Hedge}
  \begin{algorithmic}[1]
\STATE Choose positive weights $\{\lambda_t\}_{t\geq 0}$
and a non-decreasing sequence $\{\beta_t\}_{t\geq 0}$ of positive projection
parameters.
\STATE Let $a>0$, $b\in\R$, and $w_0=(1/n,\ldots,1/n)\in\R^n$.
\FOR{$t=0,\ldots,T-1$}
\STATE Set $x_t=w_t/\sum_{i=1}^nw_{t,i}$.\\
\STATE Observe loss vector $\ell_t$.\\
\STATE Set $w_{t+1,i}=\exp\left(-\frac{\lambda_t(a\ell_{t,i}+b)}{\beta_{t+1}}\right)w_t^{\frac{\beta_{t+1}}{\beta_t}}$ for any $i=1,\ldots,n$.\\
\ENDFOR
\end{algorithmic}
\end{algorithm}

Let us now discuss several strategies for choosing the weights $\gamma_t$, the projection parameters $\beta_t$, and the affine function $z\mapsto az+b$ in Algorithm \ref{alg:Extended_Hedge}. However, first we observe that the norm of each oracle return $a\ell_t+b$  and the prox-function $\de$ are bounded from above by the quantities $\boundReturns(a,b):=\max\left\{\left|-a\mu+b\right|,\left|a\rho+b\right|\right\}$ and by $D:=\ln(n)$, respectively.

\textbf{Original Hedge algorithm:} If $\beta_t=1$ and $\lambda_t=\ln(1/\gamma)$ for any $t=0,\ldots,T-1$ and with a fixed $\gamma\in (0,1)$, we recover the Hedge algorithm. This implies that the Hedge algorithm is Dual Averaging scheme.

\textbf{Optimal Hedge algorithm:} Theorem \ref{thm:ConvergenceDAS} yields for these weights and projection parameters:
\[\frac{1}{T}\max_{x\in \Delta_n}\sum_{t=0}^{T-1}\la \ell_t,x_t-x\ra\leq \frac{1}{aT\ln(1/\gamma)}\left(D+\frac{1}{2}\sum_{t=0}^{T-1}\ln^2(1/\gamma)\boundReturns^2(a,b)\right).\]
This bound is minimized by parameters $\gamma^*$, $a^*$, and $b^*$ that satisfy
\[\gamma^*=\exp\left(-\frac{2}{a^*(\mu+\rho)}\sqrt{\frac{2\ln(n)}{T}}\right)\qquad \text{and}\qquad b^*=\frac{\mu-\rho}{2}a^*,\]
where $a^*>0$. We refer to Algorithm \ref{alg:Extended_Hedge} with the just specified setting as \textit{Optimal Hedge algorithm}, for which we have by the above inequality:
\[\frac{1}{T}\max_{x\in \Delta_n}\sum_{t=0}^{T-1}\la \ell_t,x_t-x\ra\leq \frac{\mu+\rho}{2}\sqrt{\frac{\ln(n)}{T}}.\]
This result improves Bound (\ref{eq:bound_hedge}) by the additive quantity $(\mu+\rho)\ln(n)/T$ and by a multiplicative factor of $2$. Note that the resulting score function $U(z)=\left(\gamma^*\right)^{a^*z+b^*}$ does not comply with Condition (\ref{eq:condition_score_function}). Therefore, neither Theorem \ref{thm:convergenceHedge} nor its extension can be used to establish the above bound.

\textbf{Optimal Time-Independent Hedge algorithm:} The update parameter $\gamma$ depends on the number of iterations $T$ in both algorithms, the Original Hedge algorithm with the score function (\ref{eq:score_Freund}) suggested by Freund and Schapire \cite{Freund_Hedge} and the Optimal Hedge algorithm. However, when investing our money at the stock market, we might not want to fix the number of times that we adapt our portfolio in advance. We thus need an update parameter that is independent of $T$. Adapting Nesterov's strategy (\ref{eq:betas}), we choose $\gamma\in (0,1)$ and set $\lambda_t:=\ln(1/\gamma)$, $\beta_0=1$, and $\beta_{t+1}=\sum_{k=0}^{t-1}1/\beta_k$ for any $t\geq 0$. Applying Theorem \ref{thm:ConvergenceDAS}, we obtain for any $T\geq 1$:
\[\frac{1}{T}\max_{x\in \Delta_n}\sum_{t=0}^{T-1}\la \ell_t,x_t-x\ra\leq \frac{\beta_T}{aT\ln(1/\gamma)}\left(D+\frac{1}{2}\ln^2(1/\gamma)\boundReturns^2(a,b)\right).\]
We minimize the right-hand side of the above inequality, that is, we choose $a^*>0$ and set
\begin{eqnarray*}%\label{eq:optimal_parameters}
\gamma^*=\exp\left(-\frac{2\sqrt{2\ln(n)}}{a^*(\mu+\rho)}\right)\qquad \text{and}\qquad b^*=\frac{\mu-\rho}{2}a^*.
\end{eqnarray*}
Exploiting Lemma 3 in \cite{nesterov:PDS}, we obtain for the resulting method, which we refer to as the \textit{Optimal Time-Independent Hedge algorithm}, the following inequalities:
\[\frac{1}{T}\max_{x\in \Delta_n}\sum_{t=0}^{T-1}\la \ell_t,x_t-x\ra\leq \left(\mu+\rho\right)\left(\frac{1}{\left(1+\sqrt{3}\right)T}+\sqrt{\frac{2}{T}}\right)\sqrt{\frac{\ln(n)}{2}}\leq 2(\mu+\rho)\sqrt{\frac{\ln(n)}{T}}\qquad \forall\ T\geq 1.\]

\textbf{Optimal Aggressive Hedge algorithm:} The later a loss appears, the more likely it is that this loss vector contains relevant information for the future development of the investment products' performances. We conclude this section by introducing an alternative version of the Hedge algorithm, where we continuously increase the weights of the loss vectors when time proceeds. For fixed $\gamma\in (0,1)$, we set $\lambda_t=\ln(1/\gamma)(t+1)^2$ and $\beta_t=t^{2.5}$ for any $t\geq 0$. Let $T>6$. Using the relations $\sum_{t=0}^{T-1}(t+1)^2=T(T+1)(2T+1)/6 > T^3/3$, $\sum_{t=0}^{T-1}(t+1)^4\leq 2T^5/7$, and Theorem \ref{thm:ConvergenceDAS}, we obtain for Algorithm \ref{alg:Extended_Hedge}:
\begin{eqnarray*}
 \frac{6}{T(T+1)(2T+1)}\max_{x\in \Delta_n}\sum_{t=0}^{T-1}(t+1)^2\la \ell_t,x_t-x\ra&<& \frac{3}{aT^3}\left(\frac{T^{2.5}D}{\ln(1/\gamma)}+\frac{1}{2}\sum_{t=0}^{T-1}\frac{(t+1)^4\ln(1/\gamma)}{T^{2.5}}\boundReturns^2(a,b)\right)\cr
&\leq &\frac{3}{a\sqrt{T}}\left(\frac{D}{\ln(1/\gamma)}+\frac{\ln(1/\gamma)\boundReturns^2(a,b)}{7}\right).
\end{eqnarray*}
The latter quantity is minimized for $a^*$, $b^*$, and $\gamma^*$ satisfying
\begin{eqnarray*}
\gamma^*=\exp\left(-\frac{2\sqrt{7\ln(n)}}{a^*(\mu+\rho)}\right)\qquad \text{and}\qquad b^*=\frac{\mu-\rho}{2}a^*
\end{eqnarray*}
with $a^*>0$. We call the resulting method \textit{Optimal Aggressive Hedge algorithm}, for which we can rewrite the above inequality as:
\begin{eqnarray*}
 \frac{6}{T(T+1)(2T+1)}\max_{x\in \Delta_n}\sum_{t=0}^{T-1}(t+1)^2\la \ell_t,x_t-x\ra<3(\mu+\rho)\sqrt{\frac{\ln(n)}{7T}}.
\end{eqnarray*}
Note that the averaged regret reflects our time-varying choice of the weights $\lambda_t$.

\section{Numerical results}
\label{sec:Num}

\begin{table}
\begin{center}
$30$ investment products ($\mu=0.5133$, $\rho=0.5175$)\\
\begin{tabular}[b]{|l|r|r|r|r|c|}
\hline Number of iterations & $7800$ & $15600$ & $23400$ & $31200$ &  w.r.t. \\
&&&&& best product\\
[0.5ex] \hline
Best investment product & $-0.0045$ &$ -0.0034$ &$-0.0081$ &$-0.0110$& $-$ \\
Original Hedge & $0.0040$ &$0.0039$ &$-0.0020$ &$-0.0047$& $42.7\%$\\
Optimal Hedge & $0.0028$ &$0.0020$ &$-0.0042$ &$-0.0075$& $68.2\%$\\
Optimal Time-Independent Hedge & $0.0010$ &$0.0011$ &$-0.0047$ &$-0.0073$& $66.4\%$\\
Optimal Aggressive Hedge & $0.0014$ &$-0.0061$ &$-0.0183$ &$-0.0252$& $229.1\%$\\
\hline
\end{tabular}\\
\end{center}
\caption{Averaged losses achieved by the best investment product, by the Original Hedge algorithm, by the Optimal Hedge algorithm, by the Optimal Time-Independent Hedge algorithm, and by the Optimal Aggressive Hedge algorithm after one, two, three, and four months of trading. In the last column, we express the final averaged loss in percentage of the final averaged loss achieved by the best investment product.} \label{table:losses}
\end{table}

\begin{figure}[t]
\begin{center}
\includegraphics[width=\linewidth]{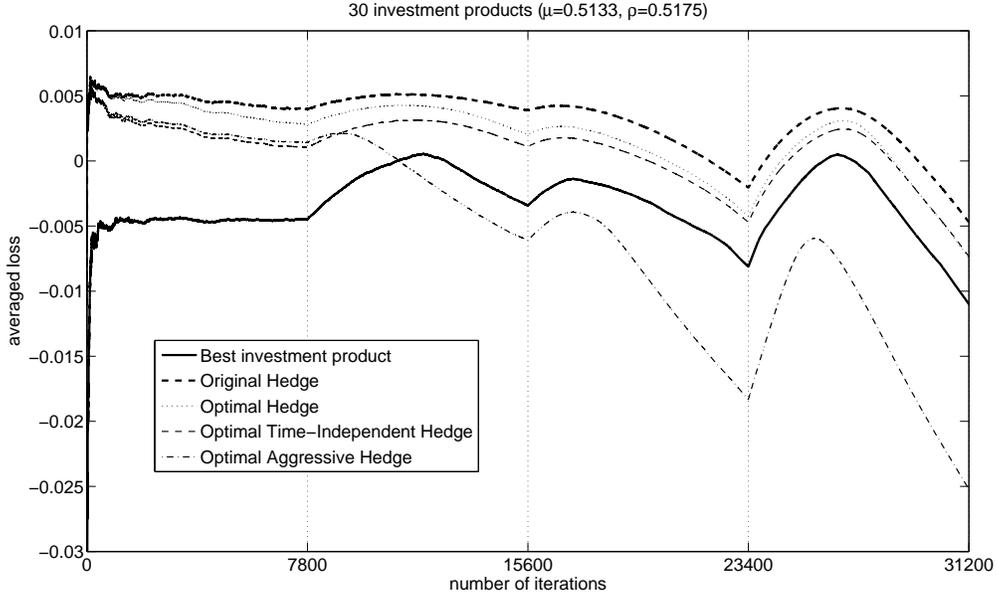}
\end{center}
\label{fig:losses} \caption{Averaged losses $\sum_{k=1}^{t}\la
\ell_{k-1},x_{k-1}\ra/t$, $t=1,\ldots,T$, achieved by the best
investment product (thick line), by the Original Hedge algorithm
(thick dashed line), by the Optimal Hedge algorithm (dotted line), by
the Optimal Time-Independent Hedge algorithm (thin dashed line), and
by the Optimal Aggressive Hedge algorithm (dashed-dotted line).}
\end{figure}

We select a pool of $n=30$ investment products and consider $T=31200$
iterations of the methods that we presented. The number $T$ is chosen
in such a way that it corresponds to the number of transactions at a
stock exchange during four months (20 trading days of 6h30 for one
month), provided that there is transaction every minute. The losses
$\ell_t\in\R^{n}$, $t=0,\ldots,T-1$, are randomly generated. The
first $7800$ losses $\left(\ell_t\right)_{t=0}^{7799}$, that is, the
losses observed during the first month, are realizations of a
multivariate normally distributed random vector with mean
$\bar{\mu}_1$ and covariance matrix $\Sigma$. The data
$(\bar{\mu}_1,\Sigma)$ is taken from \cite{database:mvn}. The losses
$\left(\ell_{7800(j-1)+k}\right)_{k=0}^{7799}$ observed in month $j$,
where $j=2,3,4$, are realizations of a multivariate normally
distributed random vector with the same covariance matrix $\Sigma$,
but with a different mean $\bar{\mu}_j$. In our experiments, we
modify each component $\bar{\mu}_{j-1,i}$ of $\bar{\mu}_{j-1}$ as
$\bar{\mu}_{j,i} = a_{j,i}\bar{\mu}_{j-1,i}+b_j$, with $b_j$ small.
The coefficient $a_{j,i}$ is negative with an increasing probability
as $j$ increases (namely $1/2$, $3/4$, and $1$), reverting the
performance of more and more products. The level of perturbation
$|a_{j,i}|$ is also increasing as $j$ increases. The experiments are
run $10$ times, and the obtained losses are averaged afterwards.

In Figure 1, we show the averaged losses, that is, $\sum_{k=1}^{t}\la
\ell_{k-1}, x_{k-1}\ra/t$ for any $t\geq 1$, achieved by the most
successful investment product at instant $t$ (obviously, this winning
product might change over time), by the Original Hedge algorithm
(with Freund and Schapire's score function as described in
(\ref{eq:score_Freund})), by the Optimal Hedge algorithm, by the
Optimal Time-Independent Hedge algorithm, and by the Optimal
Aggressive Hedge algorithm. Note that we show for the Optimal
Aggressive Hedge algorithm also the quantity $\sum_{k=1}^{t}\la
\ell_{k-1}, x_{k-1}\ra/t$, although we use a different weighting in
the algorithm and in its theoretical analysis; compare with the last
section. In Table \ref{table:losses}, we give the averaged losses
after each month.

We observe that all the extensions of the Hedge algorithm that we
suggested in this paper significantly outperform its original
counterpart. Even more interestingly, the Optimal Aggressive Hedge
algorithm achieves an averaged loss that is more than two times
better than the averaged loss of the best investment product after
$4$ months. The Optimal Aggressive Hedge algorithm outperforms the
most successful investment product, as the investment product with
the best performance has accumulated a significant loss in an early
month. This happens as we switch signs when we perturb the means of
the distribution that we use to generate random losses.

Compared to the other versions of the Hedge algorithm that we
suggested in this paper, the Optimal Hedge algorithm reacts faster
and thus more successful to the perturbations. This is due to the
increasing weights $\lambda_t$, which makes losses the more relevant
the later they appear. Recall that all the other methods consider the
losses as equally important.

\bibliographystyle{amsalpha}
\bibliography{Jordan}

\end{document}